\documentclass[11pt]{article}

\usepackage{amssymb,amscd,array}

\let\ssection=\section
\renewcommand{\section}{\setcounter{equation}{0}\ssection}

\def\d{\delta}
\def\G{\Gamma}
\def\g{\mathrm g}
\def\om{\omega}
\def\r{\rho}
\def\a{\alpha}
\def\b{\beta}
\def\s{\sigma}
\def\vfi{\varphi}

\def\l{\lambda}
\def\m{\mu}
\def\n{\nabla}

\def\implies{\Rightarrow}

\newcommand{\bbR}{\mathbb{R}}
\newcommand{\bbRP}{\mathbb{RP}}

\newcommand{\Diff}{\mathrm{Diff}}

\newcommand{\End}{\mathrm{End}}
\newcommand{\cF}{{\mathcal{F}}}
\newcommand{\cS}{{\mathcal{S}}}
\newcommand{\cD}{{\mathcal{D}}}

\newcommand{\Hom}{\mathrm{Hom}}

\newcommand{\Pol}{\mathrm{Pol}}
\newcommand{\PSL}{\mathrm{PSL}}
\newcommand{\SL}{\mathrm{SL}}
\newcommand{\Sl}{\mathrm{sl}}
\newcommand{\Og}{\mathrm{O}}

\newcommand{\Vect}{\mathrm{Vect}}

\newcommand{\cqfd}{\hspace*{\fill}\rule{3mm}{3mm}}

\begin{document}

\frenchspacing 

\def\d{\delta}
\def\om{\omega}
\def\r{\rho}
\def\a{\alpha}
\def\b{\beta}
\def\s{\sigma}
\def\vfi{\varphi}
\def\l{\lambda}
\def\m{\mu}
\def\implies{\Rightarrow}

\oddsidemargin .1truein
\newtheorem{thm}{Theorem}[section]
\newtheorem{lem}[thm]{Lemma}
\newtheorem{cor}[thm]{Corollary}
\newtheorem{pro}[thm]{Proposition}
\newtheorem{ex}[thm]{Example}
\newtheorem{rmk}[thm]{Remark}
\newtheorem{defi}[thm]{Definition}
\title{Conformal Schwarzian derivatives and conformally invariant quantization}
\author{Sofiane BOUARROUDJ\footnote{Research supported by the Japan Society 
for the Promotion of Science.}\\
 {\footnotesize   
Department of Mathematics, Keio University, Faculty of Science \& Technology}\\
{\footnotesize 3-14-1, Hiyoshi, Kohoku-ku, Yokohama, 223-8522, Japan. }\\
{\footnotesize  e-mail: sofbou@math.keio.ac.jp}
}
\date{}
\maketitle
\begin{abstract} Let $(M,\g)$ be a pseudo-Riemannian manifold. We propose a 
new approach for defining the conformal Schwarzian derivatives. These 
derivatives are 1-cocycles on the group of diffeomorphisms of $M$ related 
to the modules of linear differential operators. As operators, these 
derivatives do not depend on the rescaling of the metric $\g.$ In particular, 
if the manifold $(M,\g)$ is conformally flat, these derivatives vanish on the 
conformal group $\Og(p+1,q+1),$ where $\mathrm{dim} (M)=p+q.$ This work is a 
continuation of \cite{b2,bo2} where the Schwarzian derivative was defined on 
a manifold endowed with a projective connection.
\end{abstract}
\section{Introduction}
Let $S^1$ be the circle identified with the projective line $\bbRP^1.$ For 
any diffeomorphism $f$ of $S^1,$ the expression 
\begin{equation}
\label{clas}
S(f):= \frac{f'''(x)}{f'(x)}-\frac{3}{2} \left (\frac{f''(x)}{f'(x)}\right )^2,
\end{equation} 
where $x$ is an affine parameter on $S^1,$ is called Schwarzian derivative 
(see \cite{cara}).\\
The Schwarzian derivative has the following properties:

(i) It defines a 1-cocycle on the group of diffeomorphisms $\Diff(S^1)$ with 
values in differential quadratics (cf. \cite{ki, s}).

(ii) Its kernel is the group of projective transformations $\PSL_2(\bbR)$.\\
The aim of this paper is to propose a new approach for constructing 
the multi-dimensional conformal Schwarzian derivative. This approach was 
recently used in \cite{b2,bo2} to introduce the multi-dimensional ``projective'' 
Schwarzian derivative. The starting point of our approach is the relation between 
the Schwarzian derivative (\ref{clas}) and the space of Sturm-Liouville operators 
(see, e.g., \cite{wi}). The space of Sturm-Liouville operators is not isomorphic 
as a $\Diff(S^1)$-module to the space of differential quadratics. 
More precisely, the space of Sturm-Liouville operators is a non-trivial 
deformation of the space of differential quadratics in the sense of 
Neijenhuis and Richardson's theory of deformation (see \cite{nr}), generated 
by the 1-cocycle (\ref{clas}) (see \cite{ga} for more details). From this 
point of view, the multi-dimensional Schwarzian derivative is 
closely related to the modules of linear differential operators. To set out 
our approach, let us introduce some notation. \\
Let $M$ be a smooth manifold. We consider the space of linear differential 
operators with arguments that are $\l$-densities on $M$ and values that are 
$\mu$-densities on $M.$ We have, therefore, a two parameter family of 
$\Diff(M)$-modules denoted by $\cD_{\l,\mu}(M).$ The corresponding space of 
symbols is the space  $\cS_{\d}(M)$ of fiberwise polynomials on $T^*M$ with 
values in $\d$-densities, where $\d=\mu-\l.$ In general, the 
space $\cD_{\l,\mu}(M)$ is not isomorphic as a $\Diff(M)$-module to the space 
$\cS_{\d}(M)$ (cf. \cite{do,lo2}). However, we are interested in the following 
two cases:\\ 
(i) If $M:=\bbR^n$ is endowed with a flat projective structure (i.e. local 
action of the group $\SL_{n+1}(\bbR)$ by linear fractional transformations) 
there exists an isomorphism between $\cD_{\l,\mu}(\bbR^n)$ and 
$\cS_{\d}(\bbR^n),$ for $\d$ generic, intertwining the action of 
$\SL_{n+1}(\bbR)$ (cf. \cite{lo2}). The multi-dimensional ``projective'' 
Schwarzian derivative was defined in \cite{b2, bo2} as an obstruction to 
extend this isomorphism to the full group $\Diff(\bbR^n).$ \\
(ii) If  $M:=\bbR^n$ is endowed with a flat conformal structure (i.e. local 
action of the conformal group $\Og (p+1,q+1),$ where $p+q=n$), 
there exists an isomorphism between $\cD_{\l,\mu}(\bbR^n)$ and 
$\cS_{\d}(\bbR^n),$ for $\d$ generic, intertwining the action of 
$\Og(p+1,q+1)$ (cf. \cite{dlo,do}). In this paper we introduce the 
multi-dimensional ``conformal'' Schwarzian derivative in this context. Recall 
that in the one-dimensional case these two notions coincide in the sense that 
the conformal Lie algebra ${\mathrm o}(2,1)$ is isomorphic to the projective 
Lie algebra $\Sl_2(\bbR).$ 
\section{Differential operators and symbols}
Let $(M,\g)$ be a pseudo-Riemannian manifold of dimension $n.$ We denote by 
$\Gamma$ the Levi-Civita connection associated with the metric $\g.$ 
\subsection{Space of linear differential operators as a module}
We denote the space of tensor densities on $M$ by $\cF_{\l}(M),$ or $\cF_{\l}$ 
for simplify. This space is nothing but the space of sections of the 
line bundle $(\wedge^n T^* M)^{\otimes \l}.$ One can define in a natural way a 
$\Diff(M)$-module structure on it: for $f\in \Diff(M)$ and 
$\phi \in{\cal F}_{\l},$ in a local coordinates $(x^i)$, the action is given 
by 
\begin{equation}
\label{den}
f^*\phi=\phi\circ f^{-1}\cdot ( {J_{f^{-1}}})^{\lambda},
\end{equation}
where $J_f=\left |Df/Dx \right |$ is the Jacobian of $f$.\\ 
By differentiating this action, one can obtain the action of the Lie algebra 
of vector fields $\Vect(M).$ 
\begin{ex}{\rm 
$\cF_0= C^{\infty}(M), \, \cF_1= \Omega^{n}(M)$ (space of differential 
$n$-forms).
}
\end{ex}
Let us recall the definition of a covariant derivative on densities. If 
$\phi \in \cF_\l,$ then $\nabla \phi \in \Omega^1 (M)\otimes \cF_\l$ given 
in a local coordinates by 
$$
\nabla_i \phi= \partial_i \phi-\l \Gamma_i \phi,
$$
with $\Gamma_i= \Gamma_{ti}^t.$ (Here and bellow summation is understood over 
repeated indices).   

Consider now ${\cal D}_{\l,\mu}(M),$ the space of linear differential operators 
acting on tensor densities
\begin{equation}
A:\cF_\l\to\cF_\m.\nonumber
\label{Conv}
\end{equation}
The action of $\Diff(M)$ on $\cD_{\l,\mu}(M)$ depends on the 
two parameters $\l$ and $\m$. This action is given by the equation
\begin{equation}
f_{\l,\m}(A)=f^*\circ A\circ {f^*}^{-1},
\label{Opaction}
\end{equation}
where $f^*$ is the action (\ref{den}) of $\Diff(M)$ on $\cF_\l$.

By differentiating this action, one can obtain the action of the Lie algebra 
$\Vect(M).$

The formul{\ae} (\ref{den}) and (\ref{Opaction}) do not depend on the choice 
of the system of coordinates.

Denote by ${\cal D}^2_{\l,\mu}(M)$ the space of second-order linear 
differential operators with the $\Diff(M)$-module structure given by 
(\ref{Opaction}). The space ${\cal D}_{\l,\mu}^2(M)$ is in fact a 
$\Diff(M)$-submodule of $\cD_{\l,\m}(M).$
\begin{ex}
{\rm The space of Sturm-Liouville operators 
$\frac{d^2}{dx^2}+u(x): \cF_{-1/2}\rightarrow \cF_{3/2}$ on $S^1,$   
where $u(x)\in 
\cF_{2}$ is the potential, is a submodule of 
$\cD_{-\frac{1}{2},\frac{3}{2}}^2(S^1)$ (see \cite{wi}).
}
\end{ex}
\subsection{The module of symbols}
The space of symbols, $\Pol (T^*M),$ is the space of functions on the 
cotangent bundle $T^*M$ that are polynomials on the fibers. This space is naturally 
isomorphic to the space ${\cal S}(M)$ of symmetric contravariant tensor fields 
on $M.$ In local coordinates $(x^i,\xi_i)$, one can write $P\in {\cal S}(M)$ 
in the form  
$$P=\sum_{l\geq 0}P^{i_1, \ldots,i_l}\xi_{i_1}\cdots \xi_{i_l},$$
with $P^{i_1, \ldots,i_l}(x)\in C^{\infty}(M).$\\

We define a one parameter family of $\Diff (M)-$module on the space of 
symbols by 
$$
{\cal S}_{\d}(M):={\cal S}(M)\otimes \cF_{\d}.
$$
For $f\in \Diff(M)$ and $P\in {\cal S}_{\d}(M),$ in a local 
coordinate $(x^i)$, the action is defined by
\begin{eqnarray}
\label{actsym}
f_{\d}(P)&=& f^*P\cdot (J_{f^{-1}})^{\d},
\end{eqnarray} 
where $J_f=|Df/Dx|$ is the Jacobian of $f,$ and $f^*$ is the natural action 
of $\Diff(M)$ on ${\cal S}(M).$ 

We then have a graduation of $\Diff(M)$-modules given by 
$$
{\cal S}_\d(M)=\bigoplus_{k=0}^\infty {\cal S}_\d^k(M),
$$
where $ S_\d^k(M)$ is the space of contravariant tensor fields of degree $k$ 
endowed with the $\Diff(M)$-module structure (\ref{actsym}).

We want to study the space of contravariant tensor fields of degree 
less than two, denoted by  
${\cal S}_{\d,2}(M)$ (i.e. ${\cal S}_{\d,2}(M):={\cal S}_{\d}^2(M)\oplus 
{\cal S}_{\d}^1(M)\oplus {\cal S}_{\d}^0(M)$).

\section{Conformal Schwarzian derivatives}
Let $(M,\g)$ be a pseudo-Riemannian manifold. Denote by $\Gamma$ the 
Levi-Civita connection associated with the metric $\g.$ 
\subsection{Main definition}
\label{enfant}
It is well known that the difference between two connections is a well-defined 
tensor field of type $(2,1).$ It follows therefore that the difference 
\begin{equation}
\label{ell}
\ell(f):=f^*\Gamma-\Gamma,
\end{equation}
where $f\in \Diff(M),$ is a well-defined $(2,1)$-tensor field on $M$. 

It is easy to see that the map 
$$
f\mapsto \ell(f^{-1})
$$ defines a non-trivial 1-cocycle on $\Diff(M)$ with values in the space of 
tensor fields on $M$ of type $(2,1).$ 
 
Our first main definition is the linear differential operator 
${\cal A}(f)$ acting from ${\cal S}_{\d}^2(M)$ to ${\cal S}_\d^1(M)$ defined by
\begin{equation}
{\cal A} (f)_{ij}^k:=
{f^*}^{-1} \left( \g^{sk}\,\g_{ij}\nabla_s \right) 
-\g^{sk}\,\g_{ij}\nabla_s + c\,\left(\ell(f)^k_{ij} -\frac{1}{n} 
\, \mathrm{Sym}_{i,j}\,\delta_{i}^k \,\ell(f)_{j}\right), 
\label{MultiSchwar1}
\end{equation}
where 
\begin{equation} 
\label{qua} 
c=2-\delta n,
\end{equation} 
and $\ell(f)_{ij}^k$ are the components of the tensor (\ref{ell})
\begin{thm}
\label{main}
{\rm (i)} For all $\d \not= 2/n,$ the map $f\mapsto {\cal A}(f^{-1})$ 
defines a non-trivial 1-cocycle on $\Diff(M)$ with values in 
$\cD ({\cal S}_{\d}^2(M), {\cal S}_\d^1(M)).$ \\
{\rm (ii)} The operator (\ref{MultiSchwar1}) does not depend on the rescaling 
of the metric $\g.$ In particular, if $M:=\bbR^n$ and $\g$ is the flat 
metric of signature $p-q$, this operator vanishes on the conformal group 
$\Og(p+1,q+1).$ 
\end{thm} 
{\bf Proof.} To prove (i) we have to verify the 1-cocycle condition
$$
{\cal A}(f\circ h)={h^*}^{-1} {\cal A}(f)+ {\cal A}(h), 
\quad \mbox{for all } f,h\in \Diff(M),
$$
where $h^*$ is the natural action on 
$\cD ({\cal S}^2_{\d}(M),{\cal S}_\d^1(M)).$ This condition holds 
because the first part of the operator (\ref{MultiSchwar1}) is a coboundary 
and the second part is a 1-cocycle. 

Let us proof that this 1-cocycle is not trivial for $\d\not=2/n$. Suppose 
that there is a first-order differential operator 
$A^k_{ij}=u^{sk}_{ij}\nabla_s+v^k_{ij}$ such that 
\begin{equation}
\label{cn}
{\cal A}(f)={f^*}^{-1} A-A.
\end{equation}
From (\ref{cn}), it is easy to see that ${f^*}^{-1} v^k_{ij}-v^k_{ij}=
(2-\delta n)\left(\ell(f)^k_{ij} -\frac{1}{n} 
\, \mathrm{Sym}_{i,j}\,\delta_{i}^k \,\ell(f)_{j}\right)$. The right-hand side of this 
equation depends on the second jet of the diffeomorphism $f,$ while the 
left-hand side depends on the first jet of $f,$ which is absurd. 

For $\d=2/n,$ one can easily see that the 1-cocycle (\ref{MultiSchwar1}) 
is a coboundary. 

Let us prove (ii). Consider a metric $\tilde \g= F\cdot\g,$ where $F$ is a 
non-zero positive function. Denote by $\tilde {\cal A}(f)$ the operator 
(\ref{MultiSchwar1}) written with the metric $\tilde \g.$ We have to prove 
that $\tilde {\cal A}(f)={\cal A}(f).$ The Levi-Civita connections associated 
with  the metrics $\g$ and $\tilde \g$ are related by 
\begin{equation}
\label{lien}
\tilde \Gamma^k_{ij}=\Gamma^k_{ij}+\frac{1}{2F}\left (F_i \d^k_j +F_j \d^k_i
-F_t \,\g^{tk} \g_{ij}\right),
\end{equation} 
where $F_i=\partial_i F.$ \\
We need some formul{\ae}: denote by $\ell(f)$ the tensor (\ref{ell}) written 
with the connection $\tilde \g,$ then we have 
\begin{eqnarray}
\label{bebe}
\tilde \nabla_k P^{ij}&=&\nabla_k P^{ij}+\frac{1}{2F}\left(  \mathrm{Sym}_{i,j} 
P^{mi} \left(F_m \delta^j_k-F_t \,\g^{tj}\g_{km}\right) +(2-n\delta )\,
P^{ij}F_k\right),\\
\tilde \ell(f)^k_{ij}&=&\ell(f)^k_{ij}+\frac{1}{2 F\circ f}
\left(\mathrm{Sym}_{i,j} \stackrel{*}{F_i}\delta^k_j -\stackrel{*}{F}_t\, 
\stackrel{*}{\g}^{tk}\stackrel{*}{\g}_{ij}\right)
-\frac{1}{2F}\left( \mathrm{Sym}_{i,j} \,F_i\, \delta^k_j -F_t\, 
\g^{tk}\g_{ij}\right), 
\nonumber
\end{eqnarray}
where $\stackrel{*}{F}_i={f^{*}}^{-1}F_i$ and $\stackrel{*}{\g}_{ij}=
{f^{*}}^{-1}\g_{ij}$ for all $P^{ij}\xi_i \xi_j\in {\cal S}_{\d}^2 (M).$ \\
By substituting the formul{\ae} (\ref{bebe}) into (\ref{MultiSchwar1}) we get 
$$
{\cal A} (f)_{ij}^k=\tilde {\cal A} (f)_{ij}^k+
\frac{1}{2 F} (\delta n +c -2 ) \,{\g}^{sk} \,\g_{ij}\, F_s\, 
+  
\frac{1} {2 F\circ f}(2-c -\delta  n) 
\stackrel{*}{\g}^{tk} \stackrel{*}{\g}_{ij}\, \stackrel{*}{F}_t \, \,
\cdot
$$
We see that ${\cal A} (f)=\tilde {\cal A} (f)$ if and only if 
$c=2-\delta n .$ 

Let us prove that the operator (\ref{MultiSchwar1}) vanishes on the 
conformal group $\Og (p+1,q+1)$ in the case when $M:=\bbR^n$ is endowed 
with the flat metric $\g_0:=\mathrm{diag}(1,\ldots,1,-1,\ldots,-1)$ 
whose trace is $p-q.$ Any $f\in \Og (p+1,q+1)$ satisfies 
${{f^*}}^{-1} {\g_0}=F\cdot{\g_0},$ where $F$ is a non-zero positive 
function. This relation implies 
\begin{eqnarray}
2 \ell(f)^k_{ij} +\mathrm{Sym}_{i,j}\,\ell(f)^s_{it}\,\g_0^{tk}\,
{\g_0}_{sj}&=&\frac{1}{F}\mathrm{Sym}_{i,j}\, \partial_i F\, \delta^k_j, 
\nonumber \\ 
\mathrm{Sym}_{i,j}\, \ell(f)^s_{li} \, {\g_0}_{sj}&=&\frac{\partial_l F}{F}
{\g_0}_{ij},\nonumber \\
\ell(f)_t&=& \frac{n}{2}\frac{\partial_t F}{F}\cdot  
\nonumber
\end{eqnarray}
Sibstitute these formul{\ae} into (\ref{MultiSchwar1}). Then we obtain
 by straightforward computation that ${\cal A}(f)\equiv 0.$ \\
\cqfd

Suppose now that ${\rm dim }(M)>2.$

Our second main definition is the linear differential operator 
${\cal B}(f)$ acting from ${\cal S}^2_{\d}(M)$ to ${\cal S}^0_\d(M)$ defined by
\begin{eqnarray}
{\cal B}(f)_{ij}&=&
{f^*}^{-1} \left(\g^{st}\, \g_{ij}\nabla_s\nabla_t \right ) -\g^{st}\,\g_{ij} 
\nabla_s \nabla_t
+c_1\,\left( \ell(f)_{ij}^s -\frac{1}{n}\, 
\mathrm{Sym}_{i,j}\,\delta_{i}^s \,\ell(f)_{j} \right)\nabla_s \nonumber \\
&&+\,c_2 \, \ell(f)_i\,\ell(f)_j  +c_3 \n_{s} 
\left ( \ell(f)_{ij}^s-\frac{1}{n}
\,\mathrm{Sym}_{i,j}\, \d_{i}^s \ell_{j}(f)\right) +c_4\, \ell(f)_{ij}^s\ell(f)_s 
\nonumber  \\
&& +\,c_5\, \ell(f)_{si}^u\ell(f)_{uj}^s + c_6\left( 
 {f^{-1}}^*(R\,\g_{ij})-R\,\g_{ij}\right ),
\label{MultiSchwar2}
\end{eqnarray}
where $\ell(f)$ is the tensor (\ref{ell}), $R$ is the scalar curvature of 
the metric $\g,$ and the constants $c_1,\ldots,c_6,$ are given by 
$$
\begin{array}{ll}
c_1=2+n(1-2\delta ), & \displaystyle c_2=\frac{(2+n(1-2\delta ))(\d- 1)}{n}, 
\nonumber \\[3mm]
c_3=\displaystyle \frac{(2+n(1-2\delta ))(\d n-2)}{n-2} , &c_4=
\displaystyle \frac{(2+n(1-2\delta ))(2\d-2)}{n-2}, \\[3mm]
c_5= \displaystyle \frac{(2+n(1-2\delta ))(1-\d)n}{n-2},& c_6=
\displaystyle \frac{n(\d -1)(n\d -2)}{(n-1)(n-2)}\cdot 
\nonumber 
\label{ga}
\end{array}
$$
\begin{thm}
\label{mainp}
{\rm (i)} For all $\d \not= \frac{n+2}{2n},$ the map 
$f\mapsto {\cal B}(f^{-1})$ defines a non-trivial 1-cocycle on $\Diff(M)$ 
with values in $\cD ({\cal S}^2_{\d}(M),{\cal S}^0_\d(M))$.

{\rm (ii)} The operator (\ref{MultiSchwar2}) does not depend on the rescaling of  
the metric $\g.$ In the flat case, this operator vanishes on the conformal 
group $\Og(p+1,q+1).$ 
\end{thm}

\noindent{\bf Proof.} To prove that the map $f\mapsto {\cal B}_{ij}(f^{-1})$ 
is a 1-cocycle, one has to verify the 1-cocycle condition
$$
{\cal B}(f\circ h)=
{h^*}^{-1}{\cal B}(f)+{\cal B}(h), \quad \mbox{for all }f,h\in \Diff(M),
$$
where $h^*$ is the natural action on 
$\cD ({\cal S}^2_{\d}(M),{\cal S}_\d^0(M)).$ To do this, we use the 
formul{\ae}
\begin{eqnarray}
\label{comp}
\n_i\, f^*_{\d} P^{kl}&=&f^*_{\d}\n_i P^{kl}-\mathrm{Sym}_{k,l} \left( 
\ell(f^{-1})_{it}^k\,  f^*_{\d} \, P^{tl}\right) +\d \,\ell(f^{-1})_i\, 
f^*_{\d} P^{kl},\\[2mm]
\n_u h^* \ell(f)_{ij}^k&=&h^*\n_u \ell(f)_{ij}^k - h^*\ell(f)_{ij}^t\, 
\ell(h^{-1})_{ut}^k + \mathrm{Sym}_{i,j}\left( h^*\ell(f)_{it}^k \,
\ell(h^{-1})_{ju}^t \right), \nonumber 
\end{eqnarray}
for all $f,h \in \Diff(M)$ and for all $P^{kl}\xi_k\xi_l\in \cS_{\d}^2(M).$ \\
Let us prove that this 1-cocycle is not trivial. Suppose that there exists an 
operator $B_{ij}:= u_{ij}^{st}\nabla_s \nabla_t+ v_{ij}^s\nabla_s +t_{ij}$ 
such that 
$$
{\cal B}(f)={f^*}^{-1} B -B.
$$
It is easy to see that ${f^*}^{-1} v_{ij}^s- v_{ij}^s=
(2+n(1-2\delta))\left( \ell(f)_{ij}^s -\frac{1}{n}\, 
Sym_{i,j}\,\delta_{i}^s \,\ell(f)_{j} \right).$ The right-hand side of this 
relation depends on the second jet of $f,$ while the the left-hand side 
depends on the first jet of $f,$ which is absurd.  

For $\d= \frac{n+2}{2n},$ the 1-cocycle (\ref{MultiSchwar2}) is trivial:
${\cal B}(f)_{ij}={f^*}^{-1} (\g_{ij}B) -B\,\g_{ij},$ where 
$B:= \g^{st}\nabla_s\nabla_t -\frac{1}{4}\frac{n-2}{n-1}R$ is the so-called 
Yamabe-Laplace operator (see, e.g., \cite{besse}). \\

Let us prove (ii). Consider a metric $\tilde \g:=F\cdot \g,$ where $F$ is a 
non-zero positive function. Denote by $\tilde {\cal B}(f)$ the operator 
(\ref{MultiSchwar2}) written with the metric $\tilde \g.$ We have to prove 
that $\tilde {\cal B}(f)={\cal B}(f).$ \\
The proof is similar to the proof of part (ii) of Theorem (\ref{main}), by 
means of the equation (\ref{lien}), (\ref{bebe}) and  
\begin{eqnarray}
\tilde \nabla_l \tilde \nabla_k P^{ij}&=&\nabla_l \tilde \nabla_k P^{ij}+
\frac{1}{2F} \left( (1-\d n) F_l\,\tilde\nabla_k P^{ij}- F_k\, 
\tilde \nabla_l P^{ij}+\g^{st}\, \g_{lk}\, F_s \tilde 
\nabla_t P^{ij}\right ) \nonumber \\
&& +\frac{1}{2F} \mathrm{Sym}_{i,j} \tilde \nabla_k P^{mi} \left (F_m \, \delta^j_l -
{\g}^{sj}{\g}_{ml}\, F_s\right),\nonumber\\
\tilde \nabla_l \, \tilde \ell(f)^k_{ij}&=&\nabla_l \,\tilde \ell(f)^k_{ij} - 
\frac{1}{2F}\, F_l\,\,\tilde \ell(f)^k_{ij}+\frac{1}{2F} \left( F_t \, 
\delta^k_l -\g^{sk}\, \g_{tl}\, F_s \right ) \tilde \ell(f)^t_{ij}
\nonumber \\
&&-\frac{1}{2F} \mathrm{Sym}_{i,j} \left (F_i \delta^s_l-F_m \,{\g}^{ms}
{\g}_{il}\right)\tilde \ell(f)^k_{js}\nonumber\\
\tilde R&=& \frac{1}{F} \left ( R-(n-1)\frac{1}{F} 
\left( \g^{ij}\nabla_i F_j+ (n-6)\frac{1}{4 \, F}\g^{ij}F_i 
F_j \right) \right)  
\end{eqnarray}
for all $P^{ij}\xi_i\xi_j\in {\cal S}_\d^2(M),$ 
where $\tilde \nabla,$ $\tilde \ell(f)$ and $\tilde R$ are, the 
covariant derivative, the tensor (\ref{ell}), and the scalar curvature 
associated with the metric $\tilde \g,$ respectively.  
\subsection{Cohomology of $\Vect(M)$ and Schwarzian derivatives}
We will give here the infinitesimal 1-cocycle associated with the 1-cocycles  
${\cal A}$ and ${\cal B}$. First, let us recall the notion of a Lie derivative 
of a connection. For each $X \in  \Vect(M),$ the Lie derivative 
\begin{equation}
\label{mou}
L_X \nabla :=(Y,Z)\mapsto [X,\n_YZ]-\n_{[X,Y]} Z- \n_{Y} [X,Z] 
\end{equation}
of $\n$ is a well-known symmetric $(2,1)$-tensor field. 
The map 
$$
X\mapsto L_X {\nabla}
$$ 
defines a 1-cocycle on $\Vect(M)$ with values in the space of symmetric 
$(2,1)$-tensor fields on $M.$ \\
The linear differential operator ${\mathfrak a}$ defined by
$$
{\mathfrak a}^k_{ij}(X):=L_X \left 
(\g^{sk}\, \g_{ij}\, \nabla_s\right) 
 +c \left ((L_X \nabla)^k_{ij} 
-\frac{1}{n}\mathrm{Sym}_{i,j}\delta_j^k \, (L_X \nabla)_i\right),
$$ where the constant $c$ is as in (\ref{qua}) and $L_X \n$ is the tensor 
(\ref{mou}), acts from $\cS_{\d}^2(M)$ to $\cS_{\d}^1(M).$
The linear differential operator ${\mathfrak b}$ defined by
\begin{eqnarray}
\label{chr}
{\mathfrak b}_{ij}(X)&:=&L_X \left (\g^{st}\, \g_{ij}\, \nabla_s\nabla_t \right) 
 +c_1  \left ((L_X \nabla)^k_{ij} 
-\frac{1}{n}\mathrm{Sym}_{i,j}\delta_j^k \, (L_X \nabla_i)\right)
\nabla_k\nonumber\\
&& +c_2 \nabla_k \left ( (L_X \nabla)^k_{ij} 
-\frac{1}{n}\mathrm{Sym}_{i,j}\delta_j^k \, (L_X \nabla)_i \right )
+c_6 L_X \left ( R \, \g_{ij}\right ),
\end{eqnarray}
where the constants $c_1,c_2$ and $c_6$ are as in (\ref{ga}) and $L_X(\n)$ 
is the tensor (\ref{mou}), acts from $\cS_{\d}^2(M)$ to $\cS_{\d}^0(M).$ 

The following two propositions follow by straightforward computation.
\begin{pro}
\rm{(i)} The map $X\mapsto {\mathfrak a}^k_{ij}(X)$ defines a 
1-cocycle on $\Vect(M)$ with values in $\cD (\cS_{\d}^2(M), \cS_{\d}^1(M)).$ 

\rm{(ii)} The operator $\mathfrak a$ does not depend on the rescaling of the metric. 
In the flat case, it vanishes on the Lie algebra $\mathrm{o}(p+1,q+1),$ where 
$p+q=n.$  
\end{pro}
\begin{pro}
\rm{(i)} The map $X\mapsto {\mathfrak b}_{ij}(X)$ defines a 1-cocycle on $\Vect(M)$ 
with values in $\cD (\cS_{\d}^2(M), \cS_{\d}^0(M)).$ 

\rm{(ii)} The operator $\mathfrak b$ does not depend on the rescaling of the metric. 
In the flat case, it vanishes on the Lie algebra $\mathrm{o}(p+1,q+1),$ where 
$p+q=n.$  
\end{pro}
In section (\ref{va}), we will show that the space $\cD_{\l,\mu}^2(M)$ can 
be viewed as a non-trivial deformation of the module ${\cal S}_{2,\d}(M)$ in 
the sense of Neijenhuis and Richardson's theory of deformation 
(see also \cite{do,lo1}). 
According to the theory of deformation, the problem of ``infinitesimal'' 
deformation is related to the cohomology group 
\begin{equation}
\label{hic}
\mathrm H^1 (\Vect(M), \End ({\cal S}_{2,\d}(M))\cdot
\end{equation} 
To compute the cohomology group (\ref{hic}) we restrict the 
coefficients to the space of linear differential operators on 
${\cal S}_{2,\d}(M),$ denoted by $\cD({\cal S}_{2,\d}(M)).$ 
This space is decomposed, as a $\Vect(M)$-module, into the direct sum 
$$
\cD({\cal S}_{2,\d}(M))=\bigoplus_{k,m=0}^{2}\cD({\cal S}_{\d}^k(M), 
{\cal S}_{\d}^m(M)),
$$
where $\cD({\cal S}_{\d}^k(M), {\cal S}_{\d}^m(M))\subset 
\Hom({\cal S}_{\d}^k(M), {\cal S}_{\d}^m(M)).$

The relation between the Schwarzian derivative (\ref{clas}) and the cohomology 
group above is as follows: recall that in the one dimensional case the space 
$\cS^k_\d (S^1)$ is nothing but $\cF_{\d-k}.$ In this case, the problem of 
deformation with respect to the Lie algebra $\Sl_2(\bbR)$ is related to the 
cohomology group 
\begin{equation}
\label{ser}
{\mathrm H}^1(\Vect(S^1), \Sl_2 (\bbR); \cD (\cF_{\d-k}, \cF_{\d-l})),
\end{equation}
where $k,l=0,1,2.$ The cohomology group (\ref{ser}) was calculated in 
\cite{bo1}, it is one dimension for $k=2,\, l=0,$ and zero otherwise. The 
(unique) non-trivial class, for $k=2$ and$l=0,$ can be integrated to the group 
of diffeomorphisms $\Diff(S^1);$ it is a zero-order operator given as a 
multiplication by the Schwarzian derivative (\ref{clas}) (see \cite{bo1} for 
more details). 

In the multi-dimensional case and for $\d=0,$ the first group of differential 
cohomology of $\Vect(M)$, with coefficients in the space 
$\cD ({\cal S}^k(M),{\cal S}^m(M))$ of linear differential operators from 
${\cal S}^k(M)$ to ${\cal S}^m(M)$ was calculated in \cite{lo2}. For $n\geq2$ 
the result is as follows
\begin{equation}
\label{cal}
\mathrm H^1(\Vect(M),  \cD({\cal S}^k(M),{\cal S}^m(M)))=
\left\{
\begin{array}{ll}
\bbR\oplus \mathrm{H}_{\mathrm{DR}}^1(M), & \mbox{if}\quad  k-m=0,\\
\bbR,& \mbox{if} \quad k-m=1,m\not=0,\\
\bbR,& \mbox{if} \quad k-m=2,\\
0,&\hbox{otherwise.}
\end{array}
\right.
\end{equation}
We believe, by analogy for the one-dimensional case, that the "infinitesimal" 
multi-dimensional Schwarzian derivative is a cohomology class in the 
cohomology group above for $k-m=2.$ This class is nothing but the operator 
${\mathfrak b}$ defined in (\ref{chr}). 
\subsection{Comparison with the projective case }
\label{connex}
Let $M$ be a manifold of dimension $n.$ Fix a symmetric affine connection 
$\Gamma$ on $M$ (here $\Gamma$ is any connection not necessarily a Levi-Civita 
one). Let us recall the notion of projective connection (see \cite{kn}).

A projective connection is an equivalent class of symmetric affine connections 
giving the same unparameterized geodesics. 

Following \cite{kn}, the symbol of the projective connection is given by the 
expression
\begin{equation}
\label{connection}
\Pi_{ij}^k=\G_{ij}^k-\frac{1}{n+1}\left (\delta_i^k \G_{j}+\delta_j^k \G_{i}
\right ),
\end{equation}
where $\Gamma_{ij}^k$ are the Christoffel symbols of the connection $\Gamma$ 
and $\Gamma_i=\Gamma_{ij}^j.$\\
Two affine connection $\Gamma$ and $\tilde \Gamma$ are projectively equivalent 
if the corresponding symbols (\ref{connection}) coincide.\\
A projective connection on $M$ is called \textit{flat} if in a neighborhood of 
each point there exists a local coordinate system $(x^1,\ldots,x^n)$ such that 
the symbols $\Pi_{ij}^k$ are identically zero (see \cite{kn} for a geometric
definition). Every flat projective connection defines a projective structure 
on $M$.  \\
Let $\Pi$ and $\tilde \Pi$ be two projective connections on $M.$ The 
difference $\Pi-\tilde \Pi$ is a well-defined $(2,1)$-tensor field. 
Therefore, it is clear that a projective connection on $M$ leads to the 
following 1-cocycle on $\Diff(M)$:
\begin{equation}
{\cal C}(f^{-1})=
\left(
(f^{-1})^*\Pi_{ij}^k-\Pi_{ij}^k
\right)
dx^i\otimes dx^j\otimes\frac{\partial}{\partial x^k}
\label{ellp}
\end{equation}
This formula is independent on the choice of the coordinate system. 

By definition, the tensor field (\ref{ellp}) depends only on the projective 
class of the connection $M.$ In particular if $\Pi \equiv 0,$ this tensor 
field vanishes on the projective group $\PSL_{n+1}(\bbR).$

One can define a 1-cocycle on $\Diff(M)$ with values in 
$\cD (\cS^{2}_\d(M), \cS^{1}_\d(M))$ by contracting any symmetric 
contravariant tensor field with the tensor (\ref{ellp}). Therefore, the 
operator (\ref{MultiSchwar1}) can be viewed as the conformal analogue of the 
tensor field (\ref{ellp}). In the same spirit, the operator 
(\ref{MultiSchwar2}) can be viewed as the conformal analogue of the 
``projective'' multi-dimensional Schwarzian derivative introduced in 
\cite{b2,bo2}.
\section{Relation to the modules of differential operators}
\label{SecComput}
\subsection{Conformally equivariant quantization}
\label{tokyo}
The quantization procedure explained in this paper was first introduced in 
\cite{do,lo1}. By an equivariant quantization we mean an identification 
between the space of linear differential operators and the corresponding 
space of symbols, equivariant with respect to the action of a (finite 
dimension) sub-group $G\subset \Diff(\bbR^n).$ Recall that in the 
one-dimensional case the equivariant quantization process was carried out 
for $G=\SL_2(\bbR)$ in \cite{cmz} (see also \cite{ga}).   

The following theorems are proven in \cite{do}. 
\begin{thm}[\cite{do}]  For all $\d \not=1,$ there exists an isomorphism 
$$
{\cal Q}_{\l,\mu}: {\cal S}_{\d}^1(M)\oplus {\cal S}_{\d}^0(M)  
\rightarrow \cD_{\l,\mu}^1(M),
$$
given as follows: for all $P=P^i\xi_i+P_0\in~{\cal S}_{\d}^1(M)
\oplus {\cal S}_{\d}^0(M),$  one can associate a 
linear differential operator given by 
\begin{eqnarray}
\label{yah1}
{\cal Q}_{\l,\mu} (P)&=&
P^{i}\n_i +\a \n_i P^{i}+P_0,
\end{eqnarray}
where 
$$
\alpha=\displaystyle \frac{\l}{1-\delta}\cdot 
$$ 
This map does not depend on the rescaling of the metric, 
intertwines the action of $\Diff(M).$ 
\end{thm}
\begin{thm}[\cite{do}] For $n>2$ and for all $\d \not=\frac{2}{n}, 
\frac{n+2}{2n}, \frac{n+1}{n},\frac{n+2}{n},$ there exists an isomorphism
$$
{\cal Q}_{\l,\mu}: {\cal S}_{\d}^2(M) \rightarrow \cD_{\l,\mu}^2(M),
$$
given as follows: for all $P=P^{ij}\xi_i\xi_j \in~{\cal S}_{\d}^2(M),$ one 
can associate a linear differential operator given by 
\begin{eqnarray}
{\cal Q}_{\l,\mu} (P)&=&
P^{ij}\n_i\n_j \nonumber\\[2mm]
&&+(\a_1 \n_i P^{ij}+\a_2\,\g^{ij}\, \g_{kl}\n_i P^{kl})
\n_j 
\label{Tensor1}\\[2mm]
&&+\a_3\n_i\n_j P^{ij}+\a_4\, \g^{st}\, \g_{ij}\n_s\n_t P^{ij}
+\a_5 R_{ij}P^{ij}+\a_6 R \,\g_{ij}\,P^{ij}, 
\nonumber
\end{eqnarray}
where $R_{ij}$ (resp. $R$) are the Ricci tensor components 
(resp. the scalar curvature) of the metric $\g,$ the constants 
$\a_1,\ldots,\a_6$ are given by
$$
\begin{array}{ll}
\a_1=\displaystyle 
\frac{2(n\l+1)}
{2+n(1-\d)}, &\a_2= \displaystyle 
\frac{n(\l+\mu-1)}
{(2+n(1-\delta ))
(2-n\d)},\\[3mm]
\a_3= \displaystyle 
\frac{n\l(n\lambda+ 1)}
{(1+n(1-\delta))(2+n(1-\delta ))}, & 
\a_5= \displaystyle
\frac{n^2\l(\mu-1)}
{(n-2)(1+n(1-\delta))} ,\\[3mm]
\a_4= \displaystyle 
\frac{n\l(n^2 \mu (2-\l-\mu)+2(n\l +1)^2-n(n+1))}
{(1+n(1-\delta))(2+n(1-\d))(2+n(1-2\d))(2-n\delta )}, & 
\a_6= \displaystyle
\frac{(n\d-2)}
{(n-1)(2+n(1-2\delta))}\, \a_5 \cdot\\[3mm]
 &\\
\end{array}
$$
and has the following properties:\\
(i) It does not depend on the rescaling of the metric $\g$.\\
(ii) If $M=\bbR^n$ is endowed with a flat conformal structure, this  map 
is unique, equivariant with respect to the action of the group 
$\Og (p+1,q+1) \subset \Diff(\bbR^n)$.
\end{thm}
Before to give the formula of the conformal equivariant map in the case of 
surfaces, let us recall an interesting approach for the multi-dimensional 
Schwarzian derivative for conformal mapping \cite{os} (see also \cite{ca}). 
First, recall that all surfaces are conformally flat. This means that 
every metric can be express (locally) as
$$
\g=F^{-1}\psi^* \g_0,
$$
where $\psi$ is a conformal diffeomorphism of $M$, $F$ is a non-zero 
positive function and $\g_0$ is a metric of constant curvature. The Schwarzian 
derivative of $\psi$ is defined in \cite{os} as the following tensor 
field 
\begin{equation}
\label{osg}
S(\psi )=\frac{1}{2F}\nabla dF-\frac{3}{4F^2}\,dF\otimes dF+\frac{1}{8F^2}
\g^{-1}(dF,dF)\, \g \cdot
\end{equation}

Now we are in position to give the quantization map for the case of surfaces.

For $\d \not=1,2,\frac{3}{2},\frac{5}{2},$ and for 
each $P=P^{ij}\xi_i\xi_j +P^i\xi_i +P_0\in {\cal S}_{\d,2}(M)$  
one associates a linear differential operator given by 
\begin{eqnarray}
Q(P)&=&
P^{ij}\n_i\n_j \nonumber\\[2mm]
&&+(\a_1 \n_i P^{ij}+\a_2\,\g^{ij}\, \g_{kl}\n_i P^{kl})
\n_j \label{surf}\\[2mm]
&&+\a_3\n_i\n_j P^{ij}+\a_4\, \g^{st}\, \g_{ij}\n_s\n_t P^{ij}
\nonumber\\[2mm]
&& +\frac{4\l(\mu-1)}{2\d -3}\left( S(\psi)_{ij}P^{ij}+
\frac{1}{8(\d-1)} R \,\g_{ij}\,P^{ij}\right )+P_0, 
\nonumber
\end{eqnarray}
where $S(\psi)$ is the tensor (\ref{osg}), $R$ is the scalar curvature and 
the coefficients $\a_1,\ldots,\a_4$ are given as above. 
\begin{rmk}{\rm 
The projectively equivariant quantization map was given in \cite{lo1} (see 
also \cite{b3} for the non-flat case). The multi-dimensional projective 
Schwarzian derivative is defined as an obstruction to extend this 
isomorphism to the full group $\Diff(M).$ We will show in the next section 
that the conformal Schwarzian derivatives  defined in this paper appear as 
obstructions to extend the isomorphisms (\ref{Tensor1}), (\ref{surf}) to 
the full group $\Diff(M).$
}
\end{rmk}
\subsection{Deformation of the space of symbols ${\cal S}_{2,\d}$}
\label{va}
The goal of this section is to explicate the relation between the 1-cocycles 
(\ref{MultiSchwar1}), (\ref{MultiSchwar2}) and the space of second-order 
linear differential operators ${\cal D}_{\l,\mu}^2(M).$ Since the space 
${\cal D}_{\l,\mu}^2(M)$ is a non-trivial deformation of the space of the 
corresponding space of symbols ${\cal S}_{\d,2}(M),$ where $\d=\mu-\l$, it is 
interesting to give explicitly this deformation in term of the 1-cocycles  
(\ref{MultiSchwar1}), (\ref{MultiSchwar2}). Namely, we are looking for the 
operator $\bar f_\d={\cal Q}^{-1}_{\l,\mu} \circ f_{\l,\mu} \circ 
{\cal Q}_{\l,\mu}$ such that the diagram below is commutative
\begin{equation}
\nonumber
\begin{CD}
{\cal S}_{\d,2}(M)
@> \bar f_\d >>
{\cal S}_{\d,2}(M)
\strut\\ 
@V{\cal Q}_{\l,\mu}VV
@VV{\cal Q}_{\l,\mu}V \strut\\
\cD_{\l,\mu}^2 (M) @>  f_{\l,\mu} >> 
\cD_{\l,\mu}^2 (M) \strut
\end{CD}
\label{TheDiagram}
\end{equation}
\begin{pro}
For all $\d \not=\frac{2}{n},\frac{n+2}{2n},1,\frac{n+1}{n}, 
\frac{n+2}{n},$ the deformation of the space of symbols ${\cal S}_{\d,2}(M)$ 
by the space $\cD_{\l,\mu}^2(M)$ as a $\Diff(M)$-module is given as follows: for 
all $P=P^{ij}\xi_i \xi_j +P^i\xi_i + P_0\in {\cal S}_{\d,2}(M),$ one has   
$$
\bar f_{\d}\cdot (P  )=
T^{ij}\xi_i\xi_j+T^{i}\xi_i+T^{0},
$$ 
where 
\begin{eqnarray}
T^{ij}&=&
(f_{\d}\,P)^{ij},
\nonumber\\
\label{ExplAction}
T^i&=&
(f_{\d}\,P)^i\;+\;
\displaystyle
\frac{n(\mu+\l -1)}{(2+n(1-\d ))(2-n\delta)}\,{\cal A}^i_{kl}(f^{-1})
(f_{\d}\,P)^{kl},\\[2mm]
\displaystyle
T^0&=&
(f_\d\, P)_0
\;-\;
\displaystyle
\frac{n \l(\mu -1)}{(2+n(1-2\d))(1-\d )(1+n(1-\d))}\,{\cal B}_{kl}(f^{-1})
(f_{\d}\, P)^{kl},\nonumber
\end{eqnarray}
and $f_\d$ is the action (\ref{actsym}). 
\label{MainAct}
\end{pro}
{\bf Proof.} The proof is a simple computation using (\ref{comp}) and the 
formul{\ae}
\begin{eqnarray}
\n_i\n_j \,{f^*_{\d}}^{-1} P^{kl}&=&{f^*_{\d}}^{-1} \n_i \n_j\, P^{kl} - 
\mathrm{Sym}_{l,k}\left (
{f^*_{\d}}^{-1} \n_i\, P^{tl} \;\ell(f)_{tj}^k \right) + 
{f^*_{\d}}^{-1} \n_u\, P^{kl} \; \ell(f)_{ij}^u \nonumber \\
&& +\d \left ({f^*_{\d}}^{-1} \n_i\, P^{kl} \;\ell(f)_{j}\right) 
-\mathrm{Sym}_{k,l} \left( \n_j \ell(f)^k_{it}\; {f^*_{\d}}^{-1} 
P^{tl}\right) \\
&&+\d \n_j \ell(f)_i \; {f^*_{\d}}^{-1} P^{kl} 
-\mathrm{Sym}_{k,l}\left( \ell(f)_{it}^k \n_j \, {f^*_{\d}}^{-1} 
P^{tl}\right) +\d \ell(f)_i\,\n_j {f^*_{\d}}^{-1}P^{kl}
\nonumber\\[2mm]
\n_{i} \, f^* \phi&=& f^* \n_i \phi +\l \,\ell(f^{-1})_{i}\,  f^* \phi
\nonumber\\[2mm]
\n_j \n_{i} \, f^* \phi&=& f^* \n_j \n_i \phi +\ell(f^{-1})_{ji}^t \, f^* 
\n_t \phi +\l \mathrm{Sym}_{i,j} \,  \ell(f^{-1})_j\, f^* \n_i \phi 
\nonumber\\[2mm]
&& +\left (\l \n_j\, \ell(f^{-1})_{i}
+ \l^2 \ell(f^{-1})_{j}\, \ell(f^{-1})_{i}\right) \, f^* \phi
\nonumber\\[2mm]
f^* R_{jk}-R_{jk}&=& \n_i\, \ell(f^{-1})^i_{jk}-\n_j \, \ell(f^{-1})_k 
- \ell(f^{-1})^m_{sj}\, \ell(f^{-1})^s_{km}+
\ell(f^{-1})_m \, \ell(f^{-1})^m_{jk},\nonumber
\end{eqnarray}
for all $\phi\in \cF_\l$ and for all $f\in \Diff(M),$ where $R_{ij}$ are the 
Ricci tensor components. \\
\cqfd
\subsection{The 1-cocycle ${\cal A},$ ${\cal B}$ and the conformally 
invariant quantization}
\label{new}
We will show in this section that the quantization procedure is not invariant if 
one consider two metrics conformally equivalent. The obstruction of the 
invariance is given by the 1-coycles ${\cal A}$ and ${\cal B}.$ 

Given two conformally equivalent metrics $\g,$ $\tilde \g.$  Denote by 
$\n$, $\cal A$ and $\cal B,$ the covariant derivative, the 1-cocycles
(\ref{MultiSchwar1}) and (\ref{MultiSchwar2}) written with the metric $\g,$ 
respectively. We have 
\begin{pro}
The quantization map  ${\cal Q}_{\l,\mu}^{\g}:$ 
${\cal D}^1_{\l,\mu}(M)\rightarrow {\cal S}^1_\d(M)\oplus {\cal S}^0_\d (M)$ 
defined in (\ref{yah1}) depend only on the conformal class of the metric $\g.$
\end{pro}
{\bf Proof.} Let $\tilde \g$ be another metric conformally equivalent to 
$\g. $ That means that there exists a diffeomorphism 
$\psi : (M,\tilde \g )\rightarrow (M,\g)$ and a non-zero positive function $F$ 
such that (locally)
$$
\tilde \g=F^{-1}\cdot \psi^* \g.
$$
The Levi-Civita of the two connections are related by
\begin{equation}
{\Gamma^k_{ij}}^{ \g}= {\Gamma^k_{ij}}^{\tilde \g}+
\frac{1}{2 F} \left (\partial_i F\, 
\d_j^k 
+\partial_j F\, \d_i^k-\partial_s F\, \tilde \g^{sk}\tilde \g_{ij}\right ) -
\ell(\psi^{-1})^k_{ij},
\end{equation} 
where $\partial_i F=F_i$ and $\ell(\psi^{-1})^k_{ij}$ are the components of 
the tensor (\ref{ell}). This equation implies that 
$$
\begin{array}{lcl}
\n^{\g}_i\phi&=&\n^{\tilde \g}_i\phi-\displaystyle \frac{\l n}{2} 
\displaystyle \frac{F_i}{F}+
\l \,\ell(\psi^{-1})_i \, \phi,\\[2mm]
\n^{\g}_i P^i &=&\n^{\tilde \g}_i P^i+\displaystyle 
\frac{n(1-\d)}{2}\frac{F_i P^i}{F}-(1-\d) \ell(\psi^{-1})_{i}P^i,
\end{array}
$$
for all $\phi\in \cF_\l$ and for all $P^i\xi_i\in \cS_{\d}^1(M).$ \\
Substitute these formul{\ae} into (\ref{yah1}) we see that 
${\cal Q}^{\tilde \g}_{\l,\mu}={\cal Q}^{\g}_{\l,\mu}.$ 
\begin{pro}
The quantization map ${\cal Q}_{\l,\mu}^{\g}:$ 
${\cal D}^2_{\l,\mu}(M)\rightarrow {\cal S}^2_\d(M)$ defined in 
(\ref{Tensor1}) has the property 
$$
{\cal Q}_{\l,\mu}^{\tilde \g} (P)={\cal Q}_{\l,\mu}^{\g} (P)+d_1 \, 
{\cal A}^{s}(\psi^{-1})
(P )\nabla_s + d_2 \,\nabla_s \left ( {\cal A}^{s}(\psi^{-1})
(P) \right ) +d_3 \, {\cal B}(\psi^{-1})(P),
$$
for all $P\in {\cal S}_{\d}^2(M),$ where the constants $d_1,d_2$ and 
$d_3$ are given by 
\begin{eqnarray}
d_1&=&\frac{n(\l+\mu-1)}{(2+n(1-\d))(2-n\d))},  \quad 
\displaystyle d_2=\frac{n\l(\l+\mu-1)}{(2+n(1-\d))(2-n\d)(1-\d)}, \nonumber 
\\[1mm]
d_3&=&\displaystyle \frac{n \l(\mu -1)}{(2+n(1-2\d))(\d -1)(1+n(1-\d))} \cdot 
\nonumber 
\end{eqnarray}
\end{pro}
{\bf Proof.} The proof involves the calculation of 
$\nabla_i ^{\tilde \g}\nabla_j ^{\tilde \g} \phi,$ 
$\nabla_i ^{\tilde \g}\nabla_j ^{\tilde \g} P^{kl}$ and $R^{\tilde \g}$ 
wich is straightforward but quite complicated.  
\begin{rmk}{\rm 
The system $d_1=d_2=d_3=0$ admits a unique solution: $(\l,\mu)=(0,1).$ The 
value of $\d=1$ is called ``resonant''. In this case, the quantization map 
is not unique; there exists a one-parameter family of such isomorphism 
(see \cite{do} for more dtails.)
}
\end{rmk}
\begin{pro}
For all $f\in \Diff(M)$ and for all conformal map $\psi: (M,\tilde \g)
\rightarrow (M,\g)$, one has 

{\rm (i)} ${\cal A}^{\tilde \g}(f)=\psi^* {\cal A}^{\g}
(\psi \circ f\circ \psi^{-1}),$

{\rm (ii)} ${\cal B}^{\tilde \g}(f)=\psi^* {\cal B}^{\g} 
(\psi \circ f\circ \psi^{-1}).$
\end{pro}
{\bf Proof.} Straightforward computation.
\begin{cor}
For all conformal map $\psi (M,\tilde \g)
\rightarrow (M,\g)$ one has 

{\rm (i)} $\tilde {\cal A}^{\tilde \g}(\psi)=-{\cal A}^{\g}(\psi^{-1} ),$ 

\rm{(ii)} $\tilde {\cal B}^{\tilde \g}(\psi)=-{\cal B}^{\g}(\psi^{-1} )\cdot$ 

\end{cor}
\begin{rmk}{\rm 
The Corollary above shows that for a conformal map $\psi,$ the 1-cocycle 
$\cal B$ is still a second-order differential operator and then does not 
coincide with the Schwarzian derivative (\ref{osg}) defined by Osgood and 
Stow.
}
\end{rmk}
\section{Appendix}
We will give a formula for the Schwarzian derivative for the case of 
surfaces. As explained in section (\ref{tokyo}), all surfaces are conformally 
flat. That means that every metric can be express (locally) as
$$
\g=F^{-1}\psi^* \g_0,
$$
where $\psi$ is a conformal diffeomorphism of $M$, and $F$ is a non-zero 
positive function, $\g_0$ is a metric of constant curvature.

The explicit formula of the Schwarzian derivative in the case of 
surfaces is: the following 
\begin{eqnarray}
{\cal B}'_2(f)_{ij}&=&
{f^*}^{-1} \left( {\g}^{st}\, {\g}_{ij}\nabla_s\nabla_t \right ) 
-\g^{st}\,\g_{ij}\nabla_s \nabla_t +4(1-\delta)\left( \ell(f)_{ij}^s 
-\frac{1}{2}\,Sym_{i,j}\, 
\delta_{i}^s \,\ell(f)_{j} \right)\nabla_s \nonumber \\[2mm]
&&+4 (1-\d)^2 \ell(f)_s \left ( \ell(f)_{ij}^s -\frac{1}{4}
Sym_{i,j}\, \delta_i^s \, \ell(f)_j\right) + 
2(\d-2)(1-\d) \,Sym_{i,j}\,\n_{j}\ell_{i}(f) \nonumber\\[2mm]
&&+8\,(\d-1)^2 \left( 
 {f^{-1}}^*(S(\psi)_{ij})-S(\psi)_{ij}+\frac{1}{2}
\n_{s}\ell(f)_{ij}^s \right )+(\d-1)
 \left ({f^{-1}}^*(R\,\g_{ij})-R\,\g_{ij}\right ),
\nonumber
\label{MultiSchwar3}
\end{eqnarray}
where $S(\psi)$ is the derivative (\ref{osg}), $\ell(f)$ is the tensor 
(\ref{ell}), $R$ is the scalar curvature of the metric $\g,$ is a 
differential operator from from $\cS_{\d}^2(M)$ to $\cS_{\d}^0(M).$

Theorem (\ref{mainp}) remains true for this operator. 
\bigskip
\bigskip

{\it Acknowledgments}. It is a pleasure to acknowledge numerous fruitful 
discussions with Prof. V. Ovsienko. I am grateful to, the referee for his 
pertinent remarks, Prof. Y. Maeda and Keio University for their 
hospitality. 

\vskip 1cm


\end{document}